# Edge Addition and the Change in Kemeny's Constant

Stephen Kirkland[1], Yuqiao Li[2], John S. McAlister[3], and Xiaohong Zhang[4]

[1]Department of Mathematics, University of Manitoba, Winnipeg, MB, Canada, stephen.kirkland@umanitoba.ca.
[2]Department of Mathematics, University of Washington, Seattle, WA, USA, yuqiaoli@uw.edu.
[3]Department of Mathematics, University of Tennessee Knoxville, Knoxville, TN, USA, jmcalis6@vols.utk.edu.
[4]Centre de Recherches Mathématiques, Université de Montréal, Montréal, QC, Canada, xiaohong.zhang@umontreal.ca.

January 29, 2025

**Abstract**

Given a connected graph $G$, Kemeny's constant $\mathcal{K}(G)$ measures the average travel time for a random walk to reach a randomly selected vertex. It is known that when an edge is added to $G$, the value of Kemeny's constant may either decrease, increase, or stay the same. In this paper, we present a quantitative analysis of this behaviour when the initial graph is a tree with $n$ vertices. We prove that when an edge is added into a tree on $n$ vertices, the maximum possible increase in Kemeny's constant is roughly $\frac{2}{3}n$, while the maximum possible decrease is roughly $\frac{3}{16}n^2$. We also identify the trees, and the edges to be added, that correspond to the maximum increase and maximum decrease. Throughout, both matrix theoretic and graph theoretic techniques are employed.

## 1 Introduction and motivation

Markov chains have proven to be a much-studied family of stochastic processes, and are applied to model phenomena in a multitude of settings, including web search [22], molecular dynamics [12], and vehicle traffic networks [11]. A key quantity associated with a time-homogeneous, finite state, discrete time Markov chain is Kemeny's constant, which is well-defined whenever the transition matrix associated with the Markov chain is irreducible. Let $w = \begin{bmatrix} w_1 & \ldots & w_n \end{bmatrix}$ be the stationary distribution vector of an irreducible Markov chain and let $m_{j,k}$ be the mean first passage time from state $j$ to state $k$. For a fixed $j$, consider the sum $\sum_{k:k \neq j} m_{j,k} w_k$. Intriguingly, this quantity is independent of the choice of the state $j$ [16] and it is known as Kemeny's constant of the Markov chain. As the entries in the stationary distribution sum to 1, and $w_k = \frac{1}{m_{k,k}}, k = 1, \ldots, n$, it follows that Kemeny's constant can also be written as $\sum_{j=1}^n \sum_{k=1}^n w_j m_{j,k} w_k - 1$. From here, we observe that Kemeny's constant can also be interpreted in terms of the expected travel time from a randomly chosen initial state to a randomly chosen destination state, where both states are randomly chosen according to the stationary distribution. That last interpretation motivates the use of Kemeny's constant as a measure of the average travel time in a Markovian model for vehicle traffic networks [11].

Kemeny's constant can also be understood from a combinatorial viewpoint. In [19, Theorem 2.3], Kemeny's constant is shown to be expressible as a ratio of the sum of weights of certain spanning directed forests, to the sum of weights of certain spanning directed trees in the underlying directed graph of the corresponding transition matrix. The resulting expression for Kemeny's constant is reminiscent of the Markov Chain Tree Theorem (see [2], for example), which shows how the stationary vector can be expressed in terms of weights of spanning directed trees. Since both [19, Theorem 2.3] and the Markov Chain Tree Theorem can be derived as consequences of the All Minors Matrix Tree Theorem (see [8]), the resemblance between those two results is not so surprising. As an aside, we note that a similar combinatorial viewpoint on Kemeny's constant is adopted in [18] in order to explain the fact that $\sum_{k=1,\ldots,n, k \neq j} m_{j,k} w_k$ is independent of the choice of $j$.

An especially well-structured type of Markov chain is a random walk on a graph. Let $G$ be a connected undirected graph on vertices $1, \ldots, n$, without loops or multiple edges. The transition matrix for a random walk on $G$ can be written as $T = \Delta^{-1} A$, where $A$ is the $(0,1)$ adjacency matrix of $G$ (i.e. the $n \times n$ matrix with 1s in the positions corresponding to edges in $G$ and 0s elsewhere) and $\Delta$ is the diagonal matrix of vertex degrees. It is well-known that in this setting, the stationary distribution vector is proportional to the vector of vertex degrees.

For the random walk on our undirected graph $G$ on $n$ vertices, Kemeny's constant admits a rather simple combinatorial expression. Denote the number of edges in $G$ by $\ell$, the number of spanning trees by $\tau$, and the degree sequence by $d_1, \ldots, d_n$. For each $j, k = 1, \ldots, n$ with $j \neq k$, let $\sigma_{j,k}$ denote the number of spanning forests consisting of two trees, one of which contains vertex $j$ and the other of which contains vertex $k$. Set $\sigma_{j,j} = 0, j = 1, \ldots, n$, and let $\Sigma$ be the matrix given by $\Sigma = [\sigma_{j,k}]_{j,k=1}^n$. Corollary 2.4 in [19] states that $\mathcal{K}(G)$,

Kemeny's constant for the random walk on $G$, is given by

$$\mathcal{K}(G) = \frac{d^\top \Sigma d}{4\ell\tau}. \tag{1}$$

We note in passing that, in the case that $G$ is a tree, the corresponding matrix $\Sigma$ coincides with the distance matrix $D$ for that tree.

The last few years have seen a substantial growth in the interest in Kemeny's constant for random walks on undirected graphs; see, for example, [1, 4, 5, 6, 21, 23]. Particular attention has been paid to understanding Kemeny's constant for trees ([10, 14, 19, 24]), and to the behaviour of Kemeny's constant when an edge is added into a graph ([9, 13, 17, 20]). In the present paper, we exclusively consider finite undirected graphs without loops or multiple edges, and our goal is to analyse the effect on the corresponding value of Kemeny's constant when adding an edge to a tree. Throughout the paper, when we say adding an edge to a graph $G$, we mean connecting two distinct vertices that are not adjacent in $G$.

The following three examples set the stage for our investigation. Here, we use $K_{a,b}$ to denote the complete bipartite graph where one partite set has cardinality $a$ and the other partite set has cardinality $b$.

**Example 1.1.** Consider the graph $K_{1,3}$. The value of Kemeny's constant is $\frac{5}{2}$, and when an edge is added between two of the pendent vertices, the value of Kemeny's constant increases to $\frac{61}{24}$. This surprising effect is analogous to Braess' paradox in traffic networks; [19] explores this phenomenon further.

**Example 1.2.** Consider the graph $K_{2,2}$. The value of Kemeny's constant is equal to $\frac{5}{2}$, and when an edge is added, the value of Kemeny's constant decreases to $\frac{47}{20}$.

**Example 1.3.** For the graph $K_{8,8}$, the value of Kemeny's constant is equal to $\frac{29}{2}$. When an edge is added, the value of Kemeny's constant remains the same.

The preceding three examples show that the addition of an edge to a graph may have the effect of increasing Kemeny's constant, decreasing it, or leaving it unchanged. Our goal in this paper is to quantify the effect on Kemeny's constant when an edge is added to a connected graph. In order to do so, we focus on the family of trees. There are three reasons for doing so: i) the trees are the minimally connected graphs, and are hence a natural starting point for our investigation; ii) for the family of trees on $n$ vertices, the corresponding values of Kemeny's constant range between $\Theta(n)$ and $\Theta(n^2)$, so the trees admit a diversity of possible behaviours; iii) the combinatorial expression for Kemeny's constant is especially tractable for trees and unicyclic graphs.

Our results in this paper have potential applications to an intruder problem on a computer network, where the intruder performs a random walk on the network (graph). In this setting, Kemeny's constant measures how quickly the intruder can move through the network. The information on the change in Kemeny's constant when the graph is perturbed provides insights into what types of simple changes to the network one can make to slow down the intruder more effectively, given that our initial network is modelled by a tree.

The paper proceeds as follows. Section 2 lays out some of the necessary technical results; for a tree $T$ and a graph $G$ formed from $T$ by adding an edge, section 3 gives an attainable upper bound on $\mathcal{K}(G) - \mathcal{K}(T)$, while section 4 establishes a corresponding lower bound on $\mathcal{K}(G) - \mathcal{K}(T)$.

## 2 $\mathcal{K}(G) - \mathcal{K}(T)$

In this section, we present some technical results that will facilitate the proofs of our main results, Theorems 3.1 and 4.1. Throughout the paper, we adopt the following conventions and notation. First, our tree $T$ is on $n$ vertices. Connecting two non-adjacent vertices of $T$ creates a unique cycle, and we denote the length of that cycle by $c$ (necessarily $c \geq 3$). Without loss of generality, we take the edge to be added to be $v_1 v_c$, and denote the path in $T$ joining $v_1$ and $v_c$ by $v_1, v_2, \ldots, v_c$. Finally, we label some subtrees of $T$ as follows: $T_1$ is the subtree containing vertex $v_1$ that arises when the edge $v_1 v_2$ is deleted from $T$; $T_c$ is the subtree containing vertex $v_c$ that arises when the edge $v_{c-1} v_c$ is deleted from $T$; and for each $j = 2, \ldots, c-1$, $T_j$ is the subtree containing vertex $v_j$ that arises when the edges $v_{j-1} v_j$ and $v_j v_{j+1}$ are deleted from $T$. For each $j = 1, \ldots, c$, we use $n_j$ to denote the number of vertices in $T_j$. (In particular, $n = n_1 + \ldots + n_c$.) Figure 1 illustrates.

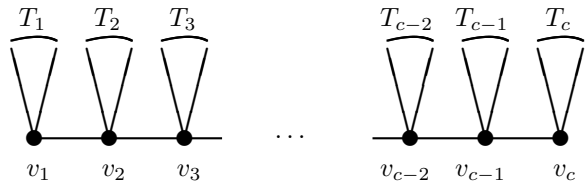

Figure 1: Our labelling convention for $T$. Each $T_j$ is on $n_j$ vertices, and includes vertex $v_j$. The edge to be added joins vertices $v_1$ and $v_c$, and the resulting graph is $G$.

Denote the distance matrix for $T$ by $D$, and the distance matrices for $T_j$ by $D_j, j = 1, \ldots, c$. Label the vertices of $T$ as follows: the vertices of $T_1$ appear first, the vertices of $T_2$ appear next, etc. Furthermore, for each $j = 1, \ldots, c$, among the vertices of $T_j$, the vertex $v_j$ appears first.

Let $\mathbf{1}$ be an all-ones vector of appropriate order, $J_{p \times q}$ be the $p \times q$ all-ones matrix, and $e_j$ be the $j$-th standard unit basis vector of appropriate order. With the above labelling, the distance matrix $D$ of $T$ can be written as a $c \times c$ block matrix, whose $(i, j)$ block is equal to

$$\begin{cases} D_i & \text{if } i = j, \\ D_i e_1 \mathbf{1}^\top + |i - j| J_{n_i \times n_j} + \mathbf{1} e_1^\top D_j & \text{otherwise.} \end{cases} \quad (2)$$

The expression (2) comes from that fact that if $i \neq j$, the distance from a vertex $a$ in $T_i$ to a vertex $b$ in $T_j$ depends only on: the distance from $a$ to $v_i$ in $T_i$ (this corresponds to the term $D_i e_1 \mathbf{1}^\top$ in (2)), the distance from $b$ to $v_j$ in $T_j$ (which corresponds to the term $\mathbf{1} e_1^\top D_j$ in (2)) and the distance from $v_i$ to $v_j$ (corresponding to the term $|i-j| J_{n_i \times n_j}$ in (2)). In matrix form,

$$D = \begin{bmatrix} D_1 - D_1 e_1 \mathbf{1}^\top - \mathbf{1} e_1^\top D_1 & & \\ & \ddots & \\ & & D_c - D_c e_1 \mathbf{1}^\top - \mathbf{1} e_1^\top D_c \end{bmatrix} \quad (3)$$

$$+ \begin{bmatrix} D_1 e_1 \\ \vdots \\ D_c e_1 \end{bmatrix} \mathbf{1}^\top + \mathbf{1} \begin{bmatrix} e_1^\top D_1 & \cdots & e_1^\top D_c \end{bmatrix} + \left[ |j-k| J_{n_j \times n_k} \right]_{j,k=1,\ldots,c}.$$

With $T$ as above we will analyse how Kemeny's constant changes when the edge $v_1 v_c$ is added to produce the resulting graph $G$. In order to do so, we first relate the distance matrix $D$ of the tree $T$ ($D$ is equal to the $\Sigma$ matrix for $T$ in (1)) to the matrix $\Sigma$ in (1) for the graph $G$.

**Lemma 2.1.** *Let $T$ be a tree as in Figure 1, and let $G$ be the graph obtained from $T$ by adding an edge between vertex $v_1$ and $v_c$. For each $a, b = 1, \ldots, n$, let $\sigma_{a,b}$ be the number of spanning forests of $G$ consisting of two trees, one of which contains vertex $a$ and the other of which contains vertex $b$. Denote the distance matrix for $T$ by $D$. The matrix $\Sigma = [\sigma_{a,b}]_{a,b=1,\ldots,n}$ satisfies $\Sigma = cD + \tilde{J}$, where $\tilde{J} = \left[ -(j-k)^2 J_{n_j \times n_k} \right]_{j,k=1,\ldots,c}.$*

*Proof.* Given distinct vertices $a, b$ of $G$, consider a spanning forest consisting of two trees, one containing $a$ and the other containing $b$. Such a forest is constructed by deleting two edges from $G$. Assume $a \in T_j$ and $b \in T_k$. Now consider the two cases: there are $(c - |j-k|)|j-k|$ such spanning forests if the two edges deleted from $G$ are both on the cycle, and $c(d_{a,b} - |j-k|)$ if one of the deleted edges is on the cycle and the other is in $T_j$ (or $T_k$ if $j \neq k$). The result follows. □

The following Lemma appears as Lemma 8.7 in [3].

**Lemma 2.2.** *Let $T$ be a tree with distance matrix $D$ and degree vector $d$. Then $d^\top D = 2 \mathbf{1}^\top D - (n-1) \mathbf{1}^\top$.*

The following is an immediate consequence of Lemma 2.2.

**Corollary 2.3.** *With $D$ and $d$ as in Lemma 2.2, we have $d^\top D d = 4 \mathbf{1}^\top D \mathbf{1} - 2(n-1)(2n-1)$.*

Now we provide an expression for the change in Kemeny's constant when an edge is added into a tree.

**Proposition 2.4.** *Let $T$ and $G$ be as in Lemma 2.1, where the subtree $T_j$ has $n_j$ vertices and distance matrix $D_j$. Let $R(T_j) = 2(n_j - 1)\mathbf{1}^\top D_j e_1 - \mathbf{1}^\top D_j \mathbf{1}$ for $j = 1, \ldots, c$. Let $B = \left[(n-1)(j-k)^2 + c|j-k|\right]_{j,k=1,\ldots,c}$ and $m = \begin{bmatrix} n_1 & \ldots & n_c \end{bmatrix}^\top$. Then*

$$\mathcal{K}(G) - \mathcal{K}(T) = \frac{1}{4cn(n-1)} \left\{ 4c \sum_{j=1}^{c} R(T_j) - 4m^\top Bm \right.$$

$$\left. + 4n(n-1)c(c-1) + 2(n-1)c^2 \right\}.$$

*Proof.* Let $D$ be the distance matrix of $T$ and $d$ be the degree vector of $T$. Since $T$ is a tree, we have $\mathcal{K}(T) = \frac{d^\top Dd}{4(n-1)}$, by (1). Let $u$ be a vector partitioned as $T$ given by $u^\top = \begin{bmatrix} e_1^\top & 0^\top & \ldots & 0^\top & e_1^\top \end{bmatrix}$, then $d + u$ is the degree vector for $G$. By equation (1) and Lemma 2.1, $\mathcal{K}(G) = \frac{(d+u)^\top \Sigma (d+u)}{4nc} = \frac{(d+u)^\top (cD + \tilde{J})(d+u)}{4nc}$, where $\tilde{J} = \left[-(j-k)^2 J_{n_j \times n_k}\right]_{j,k=1,\ldots,c}$. Consequently,

$$\mathcal{K}(G) - \mathcal{K}(T) = \frac{1}{4cn(n-1)} (-cd^\top Dd + 2(n-1)cd^\top Du + (n-1)d^\top \tilde{J} d$$

$$+ 2(n-1)d^\top \tilde{J} u + (n-1)cu^\top Du + (n-1)u^\top \tilde{J} u). \quad (4)$$

Now we evaluate the terms above. First we calculate $-d^\top Dd + 2(n-1)d^\top Du$. From Lemma 2.2 we find that $-d^\top Dd + 2(n-1)d^\top Du = 4\bigl(-\mathbf{1}^\top D\mathbf{1} + (n-1)\mathbf{1}^\top Du\bigr) + 2(n-1)$. Using (3) and simplifying now yields

$$-d^\top Dd + 2(n-1)d^\top Du = 4\left(\sum_{j=1}^{c} 2(n_j-1)\mathbf{1}^\top D_j e_1 - \mathbf{1}^\top D_j \mathbf{1}\right)$$

$$- 4\sum_{j=1}^{c}\sum_{k=1}^{c} |j-k|n_j n_k + 4n(n-1)(c-1) + 2(n-1). \quad (5)$$

Now we evaluate $d^\top \tilde{J} d$. Note that in $T$, vertex $v_j$ has one neighbour outside of $T_j$ for $j = 1, c$, and two neighbours outside $T_j$ for $j = 2, \ldots, c-1$. By the Handshaking Lemma we have

$$\sum_{k \in T_j} d_k = \begin{cases} 2n_j - 1 & \text{if } j = 1, c, \\ 2n_j & \text{otherwise.} \end{cases}$$

Therefore

$$d^\top \tilde{J} d = -\sum_{i=1}^{c}\sum_{j=1}^{c} (i-j)^2 \Bigl(\sum_{k \in T_i} d_k\Bigr)\Bigl(\sum_{\ell \in T_j} d_\ell\Bigr)$$

$$= -2(c-1)^2 + 4\sum_{j=1}^{c}((j-1)^2 + (c-j)^2)n_j - 4\sum_{j=1}^{c}\sum_{k=1}^{c}(j-k)^2 n_j n_k. \quad (6)$$

A similar computation shows

$$
\begin{aligned}
d^\top \tilde{J} u &= 2(c-1)^2 - 2\sum_{j=1}^{c}((j-1)^2 + (c-j)^2)n_j; \\
u^\top D u &= 2(c-1); \\
u^\top \tilde{J} u &= -2(c-1)^2.
\end{aligned}
\qquad (7)
$$

Now substituting equations (5), (6), and (7) into equation (4), we have

$$
\begin{aligned}
\mathcal{K}(G) - \mathcal{K}(T) &= \\
&= \frac{1}{4cn(n-1)} \Bigg\{ 4c\sum_{j=1}^{c} R(T_j) - 4c\sum_{j=1}^{c}\sum_{k=1}^{c} |j-k| n_j n_k + 4n(n-1)c(c-1) \\
&\quad + 2(n-1)c - 2(n-1)(c-1)^2 + 4(n-1)\sum_{j=1}^{c}((j-1)^2 + (c-j)^2)n_j \\
&\quad - 4(n-1)\sum_{j=1}^{c}\sum_{k=1}^{c}(j-k)^2 n_j n_k + 4(n-1)(c-1)^2 \\
&\quad - 4(n-1)\sum_{j=1}^{c}((j-1)^2 + (c-j)^2)n_j + 2(n-1)c(c-1) \\
&\quad - 2(n-1)(c-1)^2 \Bigg\} \\
&= \frac{1}{4cn(n-1)} \Bigg\{ 4c\sum_{j=1}^{c} R(T_j) - 4\sum_{j=1}^{c}\sum_{k=1}^{c}((n-1)(j-k)^2 + c|j-k|) n_j n_k \\
&\quad + 4n(n-1)c(c-1) + 2(n-1)c^2 \Bigg\}.
\end{aligned}
\qquad (8)
$$

The result follows by noting that

$$
\sum_{j=1}^{c}\sum_{k=1}^{c}((n-1)(j-k)^2 + c|j-k|) n_j n_k = m^\top B m.
$$

□

**Remark 2.5.** Fix $n_1, \ldots, n_c$ in (8), then for $j = 1, \ldots, c$, the only term that depend on the detailed structure of $T_j$ is $R(T_j) = 2(n_j - 1)\mathbf{1}^\top D_j e_1 - \mathbf{1}^\top D_j \mathbf{1}$.

Let $T$ be a tree. We finish this section with providing a combinatorial meaning of the value $R(T)$ and finding an upper and a lower bound on $R(T)$. First we study a distance-related parameter of a tree.

For two vertices $a, b$ of a tree $T$, we use $P_{a,b}$ to denote the unique path between them. The distance from a vertex $r$ to $P_{a,b}$ is the distance between vertex $r$ and the vertex $s$ on $P_{a,b}$ that is closest to $r$. When $a = b$, the path $P_{a,b}$ is the vertex $a$, and the distance from $r$ to $P_{a,b}$ is just $d_{r,a}$.

**Lemma 2.6.** *Let $T$ be a tree on $n$ vertices. For any three vertices $a, b, r$ (not necessarily distinct) of $T$, let $d_{a,b}$ denote the distance between vertices $a$ and $b$, and let $d_{r,P_{a,b}}$ denote the distance from vertex $r$ to the (unique) path $P_{a,b}$ between $a$ and $b$. Then*

$$d_{r,a} + d_{r,b} - d_{a,b} = 2d_{r,P_{a,b}}. \qquad (9)$$

*Hence*

$$0 \leq d_{r,a} + d_{r,b} - d_{a,b} = 2d_{r,P_{a,b}} \leq 2(n-1), \qquad (10)$$

*with the upper bound achieved if and only if $T$ is a path on $n$ vertices with $a = b$ being one end vertex and vertex $r$ being the other end vertex; the lower bound is achieved if and only if vertex $r$ is on the path $P_{a,b}$.*

*Proof.* Let $s$ be the vertex on $P_{a,b}$ that is closest to $r$. There are two cases to consider. If vertex $r$ is on the path $P_{a,b}$, then $r = s$ and $d_{r,P_{a,b}} = d_{r,s} = 0 = d_{r,a} + d_{r,b} - d_{a,b}$. If vertex $r$ is not on the path $P_{a,b}$, then $d_{r,P_{a,b}} = d_{r,s}$ and

$$d_{r,a} + d_{r,b} - d_{a,b} = d_{r,s} + d_{s,a} + d_{r,s} + d_{s,b} - d_{s,a} - d_{s,b} = 2d_{r,s}.$$

In either case, equation (9) holds and

$$0 \leq d_{r,P_{a,b}} = d_{r,s} \leq n - 1,$$

with lower bound achieved if and only if $r = s$, that is, $r$ is on $P_{a,b}$; with upper bound achieved if and only if $T$ is the path $P_n$ on $n$ vertices with $r$ and $s = a = b$ its two end vertices. $\square$

Now fix two vertices $r$ and $a$ in a tree $T$, and let $b$ range over all vertices of $T$, we find bounds on the sum of distances $d_{r,P_{a,b}}$, the distance from $r$ to $P_{a,b}$.

**Lemma 2.7.** *Let $T$ be a tree on $n$ vertices. For any three vertices $a, b, r$ (not necessarily distinct) of $T$, let $d_{r,P_{a,b}}$ denote the distance from vertex $r$ to the path $P_{a,b}$ between $a$ and $b$. Then for fixed vertices $r$ and $a$, we have*

$$0 \leq 2\sum_{b=1}^{n} d_{r,P_{a,b}} \leq n(n-1).$$

*Furthermore, the upper bound is achieved if and only if $T$ is a path on $n$ vertices and vertices $r$ and $a$ are the two end vertices; the lower bound is achieved if and only if $r = a$.*

*Proof.* Denote the quantity $2\sum_{b=1}^{n} d_{r,P_{a,b}}$ associated to the tree $T$ as $D(T, r, a)$, and let $d_{u,v}$ denote the distance between vertices $u$ and $v$. We first prove the upper bound by induction. If $T$ has just one vertex, then $r = a = b$, $D(T, r, a) = 0$, and the bound certainly holds. Suppose now that the inequality holds for all trees on $k$ vertices. Let $T$ be any tree on $k+1$ vertices. If $T$ has at least one pendent vertex which is neither vertex $r$ nor $a$, we take one such pendent vertex, and denote it as $k+1$ and its neighbour as $p$. Now delete the pendent vertex $k+1$ and denote the resulting tree $\tilde{T}$. Denote the distances between vertices or from a vertex to a path in $\tilde{T}$ as $\tilde{d}_{a,b}$, $\tilde{d}_{r,P_{a,b}}$, respectively. Then

$$\begin{aligned}
D(T, r, a) &= \left(2\sum_{b=1}^{k} d_{r,P_{a,b}}\right) + 2d_{r,P_{a,k+1}} \\
&= \left(2\sum_{b=1}^{k} \tilde{d}_{r,P_{a,b}}\right) + \left((\tilde{d}_{r,a}) + (\tilde{d}_{r,p}+1) - (\tilde{d}_{a,p}+1)\right) \quad \text{(by (9))} \\
&= D(\tilde{T}, r, a) + 2\tilde{d}_{r,P_{a,p}} \\
&\leq k(k-1) + 2(k-1) \quad \text{(by induction and (10))} \\
&= k^2 + k - 2 < (k+1)k.
\end{aligned}$$

Therefore if $T$ is not a path, or if $T$ is a path but at most one of vertices $r$ and $a$ is an end vertex, then the upper bound is strict. Direct computation shows that if $T$ is a path and vertices $r$ and $a$ are its two end vertices, then the upper bound is achieved.

As a sum of distances in $T$, the quantity $D(T, r, a) = 2\sum_{b \in V(T)} d_{r,P_{a,b}} \geq 0$, and it is equal to 0 if and only if $d_{r,P_{a,b}} = 0$ for all $b \in V(T)$. For a given vertex $b$, by Lemma 2.6, $d_{r,P_{a,b}} = 0$ if and only if vertex $r$ is on the path $P_{a,b}$. The condition that $r$ is on $P_{a,b}$ holds for all $b \neq r, a$ if and only if $a$ is a pendent vertex adjacent to $r$ or $a = r$. If $a \neq r$, then $d_{r,P_{a,a}} = d_{r,a} > 0$. That is, $d_{r,P_{a,b}} > 0$ for $b = a$, a contradiction. Therefore the lower bound is achieved if and only if $r = a$. $\square$

Now we provide a combinatorial meaning to $R(T)$ (used in Proposition 2.4), tight bounds for $R(T)$ and a characterization of when the bounds are achieved.

**Proposition 2.8.** *Let $T$ be a tree on $n$ vertices with distance matrix $D$. Let 1 be a fixed vertex in $T$ and let $R(T) = 2(n-1)\mathbf{1}^\top D e_1 - \mathbf{1}^\top D \mathbf{1}$. Then*

$$R(T) = 4 \sum_{1 \leq a < b \leq n} d_{1,P_{a,b}}, \qquad (11)$$

*and*

$$0 \leq R(T) \leq \frac{2}{3}n(n-1)(n-2). \qquad (12)$$

*Furthermore, the upper bound is achieved if and only if $T = P_n$ and vertex 1 is one of its end vertices; the lower bound is achieved if and only if $T$ is a star and vertex 1 is its centre.*

*Proof.* We first prove equation (11). By (9),

$$\sum_{a=1}^{n}\sum_{b=1}^{n} 2d_{1,P_{a,b}} = \sum_{a=1}^{n}\sum_{b=1}^{n}(d_{1,a} + d_{1,b} - d_{a,b}) = \sum_{a=1}^{n} 2nd_{1,a} - \sum_{a=1}^{n}\sum_{b=1}^{n} d_{a,b}$$
$$= 2n\mathbf{1}^\top D e_1 - \mathbf{1}^\top D \mathbf{1}.$$

Therefore
$$R(T) = 2(n-1)\mathbf{1}^\top D e_1 - \mathbf{1}^\top D \mathbf{1}$$
$$= \sum_{a=1}^{n}\sum_{b=1}^{n} 2d_{1,P_{a,b}} - 2\sum_{a=1}^{n} d_{1,a}$$
$$= 2 \sum_{\substack{a,b=1,\ldots,n \\ a \neq b}} d_{1,P_{a,b}}$$
$$= 4 \sum_{1 \leq a < b \leq n} d_{1,P_{a,b}}.$$

Hence $R(T) = 0$ if and only if $d_{1,P_{a,b}} = 0$ for any non-trivial path (of positive length) $P_{a,b}$ in $T$, which holds if and only if the path between any two distinct vertices passes through vertex 1, that is, $T$ is the star graph $K_{1,n-1}$ with vertex 1 being its centre.

Now we prove the upper bounds on $R(T)$ by induction on $n$ and note that for the cases $n = 1, 2$, the statement is readily verified. Suppose that (12) holds for all trees on $k$ vertices and take any tree $T$ on $k+1$ vertices. Observe that $T$ has at least one pendent vertex which is not 1. Take one such vertex, denote it by $k+1$ and its neighbour by $p$. Now delete vertex $k+1$ and denote the resulting tree $\tilde{T}$. Again denote the distances in $\tilde{T}$ as $\tilde{d}_{a,b}$ or $\tilde{d}_{r,P_{a,b}}$. First note that for any $a, u = 1, \ldots, k$,
$$d_{1,P_{a,u}} = \tilde{d}_{1,P_{a,u}}, \quad d_{1,P_{a,k+1}} = \tilde{d}_{1,P_{a,p}}. \tag{13}$$

Now by equation (11), for the tree $T$ on $k+1$ vertices,
$$R(T) = 4 \sum_{1 \leq a < b \leq k+1} d_{1,P_{a,b}}$$
$$= 4 \sum_{1 \leq a < b \leq k} d_{1,P_{a,b}} + 4 \sum_{a=1}^{k} d_{1,P_{a,k+1}}$$
$$= R(\tilde{T}) + 4 \sum_{a=1}^{k} \tilde{d}_{1,P_{a,p}} \quad \text{(by (13))}.$$

Applying the induction hypothesis on $\tilde{T}$, we have $R(\tilde{T}) \leq \frac{2}{3}k(k-1)(k-2)$, and by Lemma 2.7 we have $4\sum_{a=1}^{k} \tilde{d}_{1,P_{a,p}} \leq 2k(k-1)$. Hence
$$R(T) = R(\tilde{T}) + 4 \sum_{a=1}^{k} \tilde{d}_{1,P_{a,p}}$$
$$\leq \frac{2}{3}k(k-1)(k-2) + 2k(k-1) = \frac{2}{3}(k+1)k(k-1). \tag{14}$$

Note that equality holds in (14) if and only if $R(\tilde{T}) = \frac{2}{3}k(k-1)(k-2)$ and $4\sum_{a=1}^{k} \tilde{d}_{1,P_{a,p}} = 2k(k-1)$. By the induction hypothesis $R(\tilde{T}) = \frac{2}{3}k(k-1)(k-2)$ if and only if $\tilde{T} = P_k$ and vertex 1 is one of its end vertices. Under this condition, by Lemma 2.7, $4\sum_{a=1}^{k} \tilde{d}_{1,P_{a,p}} = 2k(k-1)$ if and only if $p$ is the other end vertex of $\tilde{T} = P_k$. Since in $T$, vertex $k+1$ is a pendent vertex adjacent to vertex $p$, the result follows. □

## 3 Upper bound on $\mathcal{K}(G) - \mathcal{K}(T)$

Let $G$ be a graph obtained from adding an edge to a tree $T$ on $n$ vertices. Assume the unique cycle in $G$ is of length $c$ (hence $c \geq 3$). In this section, we establish an attainable upper bound on $\mathcal{K}(G) - \mathcal{K}(T)$ for each fixed $c$, provide a characterization of the graphs that achieve this upper bound, and find the unique tree on $n$ vertices to which the addition of an edge provides the largest increase in $\mathcal{K}$ by analysing the dependence of the upper bound on $c$.

### 3.1 Attainable upper bound for fixed $c$

**Theorem 3.1.** *Let $T$ and $G$ be as in Lemma 2.1. Suppose that $n \geq 4$. Then*

$$\mathcal{K}(G) - \mathcal{K}(T) \leq$$
$$\frac{1}{4cn(n-1)} \left\{ 4c\frac{2}{3}(n-c+1)(n-c)(n-c-1) \right.$$
$$-8(n-c)\left(-nc\left\lfloor\frac{c+1}{2}\right\rfloor\left\lceil\frac{c+1}{2}\right\rceil + \frac{c(c+1)}{6}((2n+1)c+n-1)\right)$$
$$\left. -\frac{4(n+1)c^2(c^2-1)}{6} + 4n(n-1)c(c-1) + 2(n-1)c^2 \right\}. \tag{15}$$

*Moreover, equality holds if and only if $T$ is the tree obtained by identifying an end vertex of the path $P_{n-c+1}$ with the middle vertex (or one of the two middle vertices if $c$ is even) of $P_c$.*

*Proof.* By Proposition 2.4,

$$\mathcal{K}(G) - \mathcal{K}(T) = \frac{1}{4cn(n-1)} \left\{ 4c\sum_{j=1}^{c} R(T_j) - 4m^\top Bm \right.$$
$$\left. + 4n(n-1)c(c-1) + 2(n-1)c^2 \right\}, \tag{16}$$

where $m = \begin{bmatrix} n_1 & \dots & n_c \end{bmatrix}^\top$ and $B = \bigl[(n-1)(j-k)^2 + c|j-k|\bigr]_{j,k=1,\dots,c}$. We estimate the first two terms in the expression of $\mathcal{K}(G) - \mathcal{K}(T)$. First we find a lower bound on $m^\top Bm$. Note that

$$m^\top Bm = (m-\mathbf{1})^\top B(m-\mathbf{1}) + 2(m-\mathbf{1})^\top B\mathbf{1} + \mathbf{1}^\top B\mathbf{1}. \tag{17}$$

For each $j = 1, \ldots, c$ we have

$$\begin{aligned}
e_j^\top B\mathbf{1} &= (n-1)\sum_{k=1}^c (j-k)^2 + c\sum_{k=1}^c |j-k| \\
&= (n-1)\left(cj^2 - 2j\frac{c(c+1)}{2} + \frac{c(c+1)(2c+1)}{6}\right) \\
&\quad + c\frac{(j-1)j + (c-j)(c-j+1)}{2} \\
&= ncj(j-(c+1)) + \frac{c(c+1)}{6}((2n+1)c + n - 1). \quad (18)
\end{aligned}$$

We deduce that $e_j^\top B\mathbf{1}$ is minimized at $j = \lfloor\frac{c+1}{2}\rfloor, \lceil\frac{c+1}{2}\rceil$. Set

$$j_0 = \lfloor\frac{c+1}{2}\rfloor.$$

Since $B$ and $m - \mathbf{1}$ are both nonnegative, equation (17) implies that

$$\begin{aligned}
m^\top Bm &\geq 2(m-\mathbf{1})^\top B\mathbf{1} + \mathbf{1}^\top B\mathbf{1} \quad (19) \\
&\geq 2(n-c)e_{j_0}^\top B\mathbf{1} + \mathbf{1}^\top B\mathbf{1}, \quad (20)
\end{aligned}$$

where for the second inequality we used the fact that $e_j^\top B\mathbf{1}$ has minimum value at $j = j_0$, and that $\mathbf{1}^\top m = n$.

Since $B$ only has zero entries on the diagonal, equality holds in (19) if and only if $m - \mathbf{1} = ae_k$ for some non-negative integer $a$ and some $k \in \{1, \ldots, c\}$. From $\mathbf{1}^\top m = n$, we find that $a = n - c$. That is, for some $k$

$$m = \mathbf{1} + (n-c)e_k. \quad (21)$$

Under this condition, equality holds in (20) if and only if the index $k$ in (21) satisfies

$$k = \begin{cases} j_0 & \text{if } c \text{ is odd,} \\ j_0 \text{ or } j_0 + 1 & \text{if } c \text{ is even.} \end{cases} \quad (22)$$

Therefore $m^\top Bm$ achieves its minimum value if and only if $m = \mathbf{1} + (n-c)e_k$ with $k$ as above. Direct computation shows that the second term in (20) satisfies

$$\begin{aligned}
\mathbf{1}^\top B\mathbf{1} &= \sum_{j=1}^c \left(ncj(j-(c+1)) + \frac{c(c+1)}{6}((2n+1)c + n - 1)\right) \\
&= \frac{(n+1)c^2(c^2-1)}{6}.
\end{aligned}$$

Combing this with (18) and (20), we have

$$m^\top Bm \geq 2(n-c)\left(-nc\lfloor\tfrac{c+1}{2}\rfloor\lceil\tfrac{c+1}{2}\rceil + \tfrac{c(c+1)}{6}((2n+1)c + n - 1)\right)$$
$$+ \tfrac{(n+1)c^2(c^2-1)}{6}, \quad (23)$$

with equality holds if and only if $m = \mathbf{1} + (n-c)e_k$ for $k$ as in (22).

Now consider the term $\sum_{j=1}^{c} R(T_j)$ in $\mathcal{K}(G) - \mathcal{K}(T)$. By Proposition 2.8 and Karamata's inequality [15],

$$\sum_{j=1}^{c} R(T_j) \leq \sum_{j=1}^{c} \frac{2}{3} n_j(n_j - 1)(n_j - 2)$$

$$\leq \frac{2}{3}(n - c + 1)(n - c)(n - c - 1), \qquad (24)$$

with upper bound achieved if and only if all $n_j = 1$ apart from one, and the only non-trivial $T_j$ is a path (on $n - c + 1$ vertices) with $v_j$ being one of its end vertices.

Putting the bounds (23) and (24) in (16), we obtain the desired bound on $\mathcal{K}(G) - \mathcal{K}(T)$. By use of the conditions on when equality holds in (23) and (24), we know that equality holds in (15) if and only if $n_{j_0} = n - c + 1$ and $T_{j_0} = P_{n-c+1}$ (so $n_j = 1, j \neq j_0$) if $c$ is odd, and either $n_{j_0} = n - c + 1$ and $T_{j_0} = P_{n-c+1}$ or $n_{j_0+1} = n - c + 1$ and $T_{j_0+1} = P_{n-c+1}$ if $c$ is even (the two trees $T$ corresponding to even $c$ are isomorphic). □

Now we provide a simpler upper bound on $\mathcal{K}(G) - \mathcal{K}(T)$ by removing the floor and ceiling functions.

**Corollary 3.2.** *Suppose $n \geq 4$. Then*

$$\mathcal{K}(G) - \mathcal{K}(T) \leq$$
$$\frac{1}{2n(n-1)}\left(\frac{4n^3 - 5n^2 + 4n}{3} - c^3 - \left(\frac{n^2 - 10n}{3}\right)c^2 - (2n^2 + n)c\right).$$

*Equality holds if and only if $c$ is odd, $n_j = 1$ when $j \neq \frac{c+1}{2}$ and $T_{\frac{c+1}{2}}$ is the path on $n - c + 1$ vertices.*

*Proof.* We consider the right side of (15); since $\lfloor \frac{c+1}{2} \rfloor \lceil \frac{c+1}{2} \rceil \leq \frac{(c+1)^2}{4}$, with equality holding if and only if $c$ is odd, it follows that

$$\mathcal{K}(G) - \mathcal{K}(T) \leq$$
$$\frac{1}{4cn(n-1)}\left\{4c\frac{2}{3}(n - c + 1)(n - c)(n - c - 1)\right.$$
$$-8(n-c)\left(-nc\frac{(c+1)^2}{4} + \frac{c(c+1)}{6}((2n+1)c + n - 1)\right)$$
$$\left.-\frac{4(n+1)c^2(c^2-1)}{6} + 4n(n-1)c(c-1) + 2(n-1)c^2\right\}$$
$$= \frac{1}{2n(n-1)}\left(\frac{4n^3 - 5n^2 + 4n}{3} - c^3 - \left(\frac{n^2 - 10n}{3}\right)c^2 - (2n^2 + n)c\right).$$
$$(25)$$

□

## 3.2 Optimization over $c$

Now we maximize $\mathcal{K}(G) - \mathcal{K}(T)$ by finding the optimizing value for $c$. This gives us the unique tree on $n$ vertices to which the addiction of an edge leads to the largest increase in $\mathcal{K}$ among all trees on $n$ vertices.

**Corollary 3.3.** *Suppose that $n \geq 4$. We have*

$$\mathcal{K}(G) - \mathcal{K}(T) \leq \frac{4n^3 - 32n^2 + 85n - 81}{6n(n-1)} (\approx \frac{2n}{3} \text{ for large } n).$$

*Moreover, equality holds if and only if $T$ is the tree formed from the path on $n-1$ vertices by adding a new pendent vertex adjacent to a next-to-pendent vertex on the path, and $G$ is formed by adding the edge between the twin pendent vertices.*

*Proof.* Consider the right side of (25), and observe that since $c$ is the length of the unique cycle in $G$, we necessarily have $c \geq 3$. The right side of (25) is readily seen to be decreasing in $c$ for $c \in [3, \infty)$, and so is maximized on that interval when $c = 3$. Consequently, we have

$$\mathcal{K}(G) - \mathcal{K}(T) \leq \frac{1}{6n(n-1)}(4n^3 - 32n^2 + 85n - 81).$$

Further, for $c = 3$, equality holds in (25) if and only if $T$ is the tree formed from $P_{n-1}$ by adding a pendent vertex at a next-to-pendent vertex, and $G$ is the graph formed from $T$ by adding the edge between the twin pendent vertices of $T$. □

By Theorem 3.1, for fixed number of vertices $n$ and fixed cycle length $c$, the tree in the shape of the letter T (identifying one end vertex of $P_{n-c+1}$ to the middle vertex of $P_c$), sees the largest increase in $\mathcal{K}$ when connecting the two end vertices of $P_c$. Corollary 3.3 tells us that among these T-shaped tree, the largest increase in $\mathcal{K}$ is achieved when $c = 3$.

Now we provide a sufficient condition for Kemeny's constant to decrease when an edge is added into a tree, i.e. the analogue of Braess' paradox cannot be realized in that case.

**Corollary 3.4.** *Let $T$ be a tree on $n \geq 4$ vertices, and let $c_0$ denote the largest positive root of the equation $-c^3 - \left(\frac{n^2-10n}{3}\right)c^2 - (2n^2+n)c + \frac{4n^3-5n^2+4n}{3} = 0$. Assume that the diameter of $T$ is at least $c_0$. Construct a graph $G$ from $T$ by adding an edge to create a cycle of length $c > c_0$. Then we must have $\mathcal{K}(G) < \mathcal{K}(T)$.*

*Proof.* Let $T$ be a tree on $n$ vertices, and let $G$ be a graph obtained from $T$ by connecting two non-adjacent vertices. Assume the unique cycle in $G$ is of length $c$. Observe from (25) that if

$$-c^3 - \left(\frac{n^2-10n}{3}\right)c^2 - (2n^2+n)c + \frac{4n^3-5n^2+4n}{3} < 0, \qquad (26)$$

then necessarily $\mathcal{K}(G) < \mathcal{K}(T)$. The result follows from that (26) holds for $c > c_0$. □

**Remark 3.5.** We investigate, for large $n$, the values of $c$ for which (26) holds. Set $g(c) = -c^3 - \left(\frac{n^2-10n}{3}\right)c^2 - (2n^2+n)c + \frac{4n^3-5n^2+4n}{3}$; by Descartes' rule of signs we know that for $n \geq 10$, there is a unique $c_0 > 0$ such that $g(c_0) = 0$. Computations reveal that $g(2\sqrt{n}) = -4n^{\frac{5}{2}} + O(n^2) < 0$ while $g(2\sqrt{n} - 3) = \frac{44}{3}n^2 + O(n^{\frac{3}{2}}) > 0$. Consequently, we find that $2\sqrt{n} - 3 < c_0 < 2\sqrt{n}$ for large values of $n$. In particular, when $n$ is large, if $2\sqrt{n} < c$, then necessarily $c > c_0$ and hence $\mathcal{K}(G) < \mathcal{K}(T)$.

## 4 Lower bound on $\mathcal{K}(G) - \mathcal{K}(T)$

This section provides an attainable lower bound on $\mathcal{K}(G) - \mathcal{K}(T)$ for each $c$ and a characterization of the graphs that achieve this upper bound, as well as an analysis of the dependence of the lower bound on $c$, from which we find the optimal $c$ to minimize this lower bound.

### 4.1 Attainable lower bound for fixed $c$

**Theorem 4.1.** *Let $T$ the $G$ be as in Lemma 2.1. Suppose that $n \geq 4$. Then*

$$\begin{aligned}
\mathcal{K}(G) - \mathcal{K}(T) \geq \frac{1}{4cn(n-1)} &\left\{ -\frac{2}{3}(n+1)c^2(c^2-1) \right. \\
-8(n-c)&\left(-nc^2 + \frac{c(c+1)}{6}((2n+1)c+n-1)\right) \\
-8\left\lfloor\frac{n-c}{2}\right\rfloor &\left\lceil\frac{n-c}{2}\right\rceil ((n-1)(c-1)^2 + c(c-1)) \\
+4n(n-1)&\left.c(c-1) + 2(n-1)c^2 \right\},
\end{aligned} \tag{27}$$

*with equality holding if and only if $T$ is formed from the path on $c$ vertices by appending $\lfloor\frac{n-c}{2}\rfloor$ pendent vertices at one end point, and $\lceil\frac{n-c}{2}\rceil$ vertices at the other end point, and $G$ is obtained from $T$ by connecting the two vertices being appended.*

*Proof.* By Proposition 2.4,

$$\mathcal{K}(G) - \mathcal{K}(T) =$$
$$\frac{1}{4cn(n-1)}\left\{4c\sum_{j=1}^{c} R(T_j) - 4m^\top Bm + 4n(n-1)c(c-1) + 2(n-1)c^2\right\},$$

where $B = \left[(n-1)(j-k)^2 + c|j-k|\right]_{j,k=1,\ldots,c}$ and $m = \begin{bmatrix} n_1 & \ldots & n_c \end{bmatrix}^\top$. We first check the term $\sum_{j=1}^{c} R(T_j)$. By Proposition 2.8, $\sum_{j=1}^{c} R(T_j) \geq 0$, with equality holds if and only if for each $j = 1, \ldots, c$, $T_j$ is a star on $n_j$ vertices with vertex $v_j$ its centre. Since the term $m^\top Bm$ only depends on the values $n_1, \ldots, n_c$, not on the structures of $T_1, \ldots, T_c$, we know that when $\mathcal{K}(G) - \mathcal{K}(T)$ is minimized, each $T_j$ is as above. Henceforth we assume that is the case. Therefore

$$\mathcal{K}(G) - \mathcal{K}(T) = \frac{1}{4cn(n-1)}\Big(-4m^\top Bm + 4n(n-1)c(c-1) + 2(n-1)c^2\Big).$$

We now prove that when $m^\top Bm$ is maximized, all $n_j = 1$ apart from $j = 1, c$. Suppose on the contrary that $n_j \geq 2$ for some $2 \leq j \leq c-1$, and consider $\tilde{m}^\top B\tilde{m}$, where $\tilde{m} = m + (e_{j-1} - e_j)$, and $\hat{m}^\top B\hat{m}$ where $\hat{m} = m + (e_{j+1} - e_j)$. We find that

$$\tilde{m}^\top B\tilde{m} = m^\top Bm + 2(e_{j-1} - e_j)^\top Bm + (e_{j-1} - e_j)^\top B(e_{j-1} - e_j),$$

and

$$\hat{m}^\top B\hat{m} = m^\top Bm + 2(e_{j+1} - e_j)^\top Bm + (e_{j+1} - e_j)^\top B(e_{j+1} - e_j).$$

We now check the three types of summands above. Observe that $(e_{j-1} - e_j)^\top B(e_{j-1} - e_j) = -2(n-1+c)$ and $(e_{j+1} - e_j)^\top B(e_{j+1} - e_j) = -2(n-1+c)$. Also,

$$(e_{j-1} - e_j)^\top Bm =$$
$$\sum_{k=1}^{c} \Big((n-1)(j-1-k)^2 + c|j-1-k| - (n-1)(j-k)^2 - c|j-k|\Big)n_k =$$
$$\sum_{k=1}^{j-1} \Big(-(n-1)(2j-2k-1) - c\Big)n_k + (n-1+c)n_j$$
$$+ \sum_{k=j+1}^{c} \Big(-(n-1)(2j-2k-1) + c\Big)n_k.$$

Similarly, we find that

$$(e_{j+1} - e_j)^\top Bm = \sum_{k=1}^{j-1} \Big((n-1)(2j-2k+1) + c\Big)n_k + (n-1+c)n_j$$
$$+ \sum_{k=j+1}^{c} \Big((n-1)(2j-2k+1) - c\Big)n_k.$$

From the above we find that

$$(\tilde{m}^\top B\tilde{m} - m^\top Bm) + (\hat{m}^\top B\hat{m} - m^\top Bm) =$$
$$2\sum_{k=1}^{j-1}(n-1)n_k + 2(n-1+c)n_j + 2\sum_{k=j+1}^{c}(n-1)n_k - 4(n-1+c) =$$
$$2((n-1)(n-2) + c(n_j - 2)) > 0.$$

We deduce that one of $\tilde{m}^\top B\tilde{m} - m^\top Bm$ and $\hat{m}^\top B\hat{m} - m^\top Bm$ is positive. Consequently, when $\mathcal{K}(G) - \mathcal{K}(T)$ is minimized, it must be the case that $n_j = 1, j = 2, \ldots, c-1$, $T_1$ is a star with centre $v_1$, and $T_c$ is a star with centre $v_c$. Hence $m = \mathbf{1} + \ell e_1 + (n - \ell - c)e_c$ for some $0 \leq \ell \leq n - c$. For this $m$, a computation shows that $m^\top B m = \frac{(n+1)c^2(c^2-1)}{6} + 2(n-c)\left(-nc^2 + \frac{c(c+1)}{6}((2n+1)c + n - 1)\right) + 2\ell(n - \ell - c)\big((n-1)(c-1)^2 + c(c-1)\big)$. Evidently this is maximized exactly when $\ell = \lfloor \frac{n-c}{2} \rfloor, \lceil \frac{n-c}{2} \rceil$. It now follows that (27) holds with equality holding if and only if $T$ is formed from $P_c$ by appending $\lfloor \frac{n-c}{2} \rfloor$ pendent vertices at one end point, and $\lceil \frac{n-c}{2} \rceil$ vertices at the other end point. □

### 4.2 Optimization over $c$

Now we minimize $\mathcal{K}(G) - \mathcal{K}(T)$ by finding the optimizing value for $c$. Unlike the upper bound, the removal of the floor and ceiling functions in the lower bound (the right side of (27)) does not give polynomials in $c$, and the analysis in this bound is more involved.

Reformulating (27), we have the following. If $n - c$ is even then

$$\mathcal{K}(G) - \mathcal{K}(T) \geq \frac{2}{12cn(n-1)}[c^4 + (2n^2 - 2n - 3)c^3 + (-3n^3 + 3n - 1)c^2$$
$$+ (6n^3 - 5n^2 + 2n)c - 3n^3 + 3n^2]. \tag{28}$$

If $n - c$ is odd then

$$\mathcal{K}(G) - \mathcal{K}(T) \geq \frac{2}{12cn(n-1)}[c^4 + (2n^2 - 2n - 3)c^3 + (-3n^3 + 6n - 1)c^2$$
$$+ (6n^3 - 5n^2 - 4n + 3)c - 3n^3 + 3n^2 + 3n - 3]. \tag{29}$$

Our goal now is to find the positive values of $c$ that minimize the right hand sides of (28) and (29), respectively. Equivalently, we see the positive values of $c$ that minimize $f_e(c)$ and $f_o(c)$ defined by

$$f_e(c) = c^3 + (2n^2 - 2n - 3)c^2 + (-3n^3 + 3n - 1)c + (6n^3 - 5n^2 + 2n)$$
$$+ \frac{-3n^3 + 3n^2}{c}, \tag{30}$$
$$f_o(c) = c^3 + (2n^2 - 2n - 3)c^2 + (-3n^3 + 6n - 1)c + (6n^3 - 5n^2 - 4n + 3)$$
$$+ \frac{-3n^3 + 3n^2 + 3n - 3}{c}. \tag{31}$$

Assume $c \neq 0$. Let

$$g_e(c) = 3c^4 + (4n^2 - 4n - 6)c^3 + (-3n^3 + 3n - 1)c^2 + 3n^3 - 3n^2,$$

that is, $g_e(c) = c^2 f'_e(c)$. Similarly, let

$$g_o(c) = 3c^4 + (4n^2 - 4n - 6)c^3 + (-3n^3 + 6n - 1)c^2 + 3n^3 - 3n^2 - 3n + 3,$$

that is, $g_o(c) = c^2 f'_o(c)$. We will find the desired minimizing values of $c$ via the following sequence of claims about the roots of $g_e(c)$ and of $g_o(c)$, in particular, about their largest (positive) roots. The first two claims provide us with some intervals where the roots of the two functions lie.

**Claim 1**: *Suppose that $n \geq 3$. Then $g_e$ has a root in each of the following intervals:* $(-\infty, -2), (-2, 0), (0, 2), (2, \infty)$.
To see this, observe that $g_e(-2) = -9n^3 - 35n^2 + 44n + 92 < 0, g_e(0) = 3n^3 - 3n^2 > 0, g_e(2) = -9n^3 + 29n^2 - 20n - 4 < 0$. The conclusion now follows from the Intermediate Value Theorem.

**Claim 2**: *Suppose that $n \geq 3$. Then $g_o$ has a root in each of the following intervals:* $(-\infty, -2), (-2, 0), (0, 2), (2, \infty)$.
We have $g_o(-2) = -9n^3 - 35n^2 + 53n + 95 < 0, g_o(0) = 3n^3 - 3n^2 - 3n + 3 > 0, g_o(2) = -9n^3 + 29n^2 - 11n - 1 < 0$. We then apply the Intermediate Value Theorem.

Observe that for $c > 0$, $f'_e(c)$ is positive, negative, or zero if and only if $g_e(c)$ is. It now follows that the root of $g_e$ in $(0, 2)$ is a local maximum for $f_e$, i.e. for the right hand side of (28). Similarly, the root of $g_o$ in $(0, 2)$ is a local maximum for the right hand side of (29). We deduce that, on the interval $(0, n)$, the largest positive root of $g_e$, respectively $g_o$, yields the global minimum for the lower bound on $\mathcal{K}(G) - \mathcal{K}(T)$. Denote the largest positive root of $g_e$ by $c_e$, and the largest positive root of $g_o$ by $c_o$. Now we find a small interval that contains these two roots. It provides a good approximation to the two roots as $n$ gets large.

**Claim 3**: *For $n \geq 3$ we have*

$$2 < \frac{3n}{4} + \frac{21}{64} - \frac{7}{4n} < c_e, c_o < \frac{3n}{4} + \frac{21}{64}.$$

We again rely on the Intermediate Value Theorem. We have the following observations.
i)
$$g_e\left(\frac{3n}{4} + \frac{21}{64}\right) = \frac{4605}{2048}n^3 - \frac{337233}{65536}n^2 - \frac{379071}{262144}n - \frac{4779117}{16777216}.$$
Considering that expression as a polynomial in $n$, we find from Descartes' rule of signs that it has at most one real positive root. Evaluating at $n = 0, n = 3$ we deduce that there is a root in $(0, 3)$. Hence $g_e\left(\frac{3n}{4} + \frac{21}{64}\right) > 0$ for $n \geq 3$.
ii)
$$g_e\left(\frac{3n}{4} + \frac{21}{64} - \frac{7}{4n}\right) =$$
$$\frac{1}{16777216n^4}(-28336128n^7 - 152391936n^6 + 389562432n^5 + 454266771n^4$$
$$- 582546496n^3 - 704847360n^2 + 185450496n + 472055808).$$

Again using Descartes' rule of signs, it follows that each of $-28336128n^7 - 152391936n^6 + 389562432n^5 + 454266771n^4$ and $-582546496n^3 - 704847360n^2 + 185450496n + 472055808$ has precisely one positive root, in $(0, 3)$. Consequently, $g_e\left(\frac{3n}{4} + \frac{21}{64} - \frac{7}{4n}\right) < 0$ for $n \geq 3$.

iii)
$$g_o\left(\frac{3n}{4} + \frac{21}{64}\right) = \frac{8061}{2048}n^3 - \frac{240465}{65536}n^2 - \frac{1080831}{262144}n + \frac{45552531}{16777216}.$$

Evaluating at $n = -1, 0, \frac{3}{5}, 2$, we find that this expression has one roots in each of $(-1, 0), (0, \frac{3}{5})$ and $(\frac{3}{5}, 2)$. Hence $g_o\left(\frac{3n}{4} + \frac{21}{64}\right) > 0$ for $n \geq 3$.

iv)
$$g_o\left(\frac{3n}{4} + \frac{21}{64} - \frac{7}{4n}\right) =$$
$$\frac{1}{16777216n^4}(-24576n^7 - 127619328n^6 + 212529216n^5 + 446795667n^4$$
$$-428405824n^3 - 704847360n^2 + 185450496n + 472055808).$$

From Descartes' rule of signs, each of $-24576n^7 - 127619328n^6 + 212529216n^5 + 446795667n^4$ and $-428405824n^3 - 704847360n^2 + 185450496n + 472055808$ has precisely one positive root, in $(0, 3)$. Consequently, $g_o\left(\frac{3n}{4} + \frac{21}{64} - \frac{7}{4n}\right) < 0$ for $n \geq 3$.

From i), ii) and the Intermediate Value Theorem, we find that $\frac{3n}{4} + \frac{21}{64} - \frac{7}{4n} < c_e < \frac{3n}{4} + \frac{21}{64}$, while iii), iv) and the Intermediate Value Theorem yield $\frac{3n}{4} + \frac{21}{64} - \frac{7}{4n} < c_o < \frac{3n}{4} + \frac{21}{64}$.

**Claim 4**: *Suppose that $n \geq 23$. Then*

$$\lfloor c_e \rfloor = \left\lfloor \frac{3n}{4} + \frac{21}{64} \right\rfloor, \lfloor c_o \rfloor = \left\lfloor \frac{3n}{4} + \frac{21}{64} \right\rfloor, \lceil c_e \rceil = \left\lceil \frac{3n}{4} + \frac{21}{64} \right\rceil, \lceil c_o \rceil = \left\lceil \frac{3n}{4} + \frac{21}{64} \right\rceil.$$

We analyse these items by cases, according to the remainder of $n$ modulo 4. Suppose that $n = 4k$ for some $k \in \mathbb{N}$. We have $\frac{3n}{4} + \frac{21}{64} = 3k + \frac{21}{64}$, so that $\lfloor \frac{3n}{4} + \frac{21}{64} \rfloor = 3k$. Also, note that $\frac{3n}{4} + \frac{21}{64} - \frac{7}{4n} = 3k + \frac{21}{64} - \frac{7}{4n} \geq 3k$, the last since $n \geq 8$. From $3k \leq \frac{3n}{4} + \frac{21}{64} - \frac{7}{4n} < c_e, c_o < \frac{3n}{4} + \frac{21}{64} < 3k+1$, we conclude that $\lfloor c_e \rfloor = \lfloor \frac{3n}{4} + \frac{21}{64} \rfloor = \lfloor c_o \rfloor$ and $\lceil c_e \rceil = \lceil \frac{3n}{4} + \frac{21}{64} \rceil = \lceil c_o \rceil$.

Next suppose that $n = 4k+1$ for some $k \in \mathbb{N}$. We have $\frac{3n}{4} + \frac{21}{64} = 3k + \frac{69}{64}$, so that $\lfloor \frac{3n}{4} + \frac{21}{64} \rfloor = 3k+1$. Also, note that $\frac{3n}{4} + \frac{69}{64} - \frac{7}{4n} = 3k+1 + \frac{5}{64} - \frac{7}{4n} \geq 3k+1$, the last since $n \geq 23$. We conclude that the claim holds. The remaining cases are established similarly.

**Claim 5**: *For $n \geq 23$, we have*

$$n - \lfloor c_e \rfloor = n - \lfloor c_o \rfloor \equiv \begin{cases} 0 \pmod{2} & \text{if } n \equiv 0, 1, 6, 7 \pmod{8}, \\ 1 \pmod{2} & \text{if } n \equiv 2, 3, 4, 5 \pmod{8}. \end{cases}$$

We analyse according to the value of $n \pmod 4$. We have

$$\lfloor c_e \rfloor = \lfloor c_o \rfloor = \begin{cases} 3k & \text{if } n = 4k, k \in \mathbb{N}, \\ 3k+1 & \text{if } n = 4k+1, k \in \mathbb{N}, \\ 3k+1 & \text{if } n = 4k+2, k \in \mathbb{N}, \\ 3k+2 & \text{if } n = 4k+3, k \in \mathbb{N}. \end{cases}$$

We also have

$$n - \lfloor c_e \rfloor = n - \lfloor c_o \rfloor = \begin{cases} k & \text{if } n = 4k, k \in \mathbb{N}, \\ k & \text{if } n = 4k+1, k \in \mathbb{N}, \\ k+1 & \text{if } n = 4k+2, k \in \mathbb{N}, \\ k+1 & \text{if } n = 4k+3, k \in \mathbb{N}. \end{cases}$$

Now the result follows.

Let $c^* = \frac{3n}{4} + \frac{21}{64}$. Note that $c^* \notin \mathbb{Z}$. Therefore $\lceil c^* \rceil \neq \lfloor c^* \rfloor$, and they are of opposite parity. From the considerations above we have the following.

**Theorem 4.2.** *Suppose that $n \geq 23$ and let $c^* = \frac{3n}{4} + \frac{21}{64}$. Let $f_e(c)$ and $f_o(c)$ be defined as in equations (30) and (31), respectively. Then for $n \equiv 0, 1, 6, 7 \pmod 8$, the minimum value of $\mathcal{K}(G) - \mathcal{K}(T)$ is*

$$\frac{1}{6n(n-1)} \min\{f_e(\lfloor c^* \rfloor), f_e(\lfloor c^* \rfloor + 2), f_o(\lceil c^* \rceil), f_o(\lceil c^* \rceil - 2)\}.$$

*For $n \equiv 2, 3, 4, 5 \pmod 8$, the minimum value of $\mathcal{K}(G) - \mathcal{K}(T)$ is*

$$\frac{1}{6n(n-1)} \min\{f_o(\lfloor c^* \rfloor), f_o(\lfloor c^* \rfloor + 2), f_e(\lceil c^* \rceil), f_e(\lceil c^* \rceil - 2)\}.$$

By Theorem 4.1, for fixed number of vertices $n$ and for fixed cycle length $c$, trees obtained from two stars of the size $\lfloor \frac{n-c}{2} \rfloor$ and $\lceil \frac{n-c}{2} \rceil$, respectively, by connecting the two centres with a path of length $c - 1$, sees the largest decrease in $\mathcal{K}$ when the two centres are connected with an edge. Theorem 4.2 provides us four possible values of $c$, which are all close to $\lfloor c^* \rfloor = \lfloor \frac{3n}{4} + \frac{21}{64} \rfloor$, for the biggest decrease in $\mathcal{K}$ to occur.

**Remark 4.3.** Since $c^* = \frac{3n}{4} + O(1)$, it follows that $\mathcal{K}(G) - \mathcal{K}(T) \geq -\frac{3n^2}{16} + O(n)$, and when $T$ is formed from $P_{\lfloor c^* \rfloor}$ by adding $\lfloor \frac{n-\lfloor c^* \rfloor}{2} \rfloor$ pendent vertices at one end vertex and $\lceil \frac{n-\lfloor c^* \rfloor}{2} \rceil$ pendent vertices at the other end vertex, $\mathcal{K}(G) - \mathcal{K}(T) = -\frac{3n^2}{16} + O(n)$. Note that for this $T$ we have $\mathcal{K}(T) = \frac{39}{128}n^2 + O(n)$, while for the graph $G$ formed by adding the edge $v_1 v_{\lfloor c^* \rfloor}$, we have $\mathcal{K}(G) = \frac{15}{128}n^2 + O(n)$. In particular, for large $n$ we see that the addition of a single edge to $T$ yields a graph $G$ for which Kemeny's constant is less than 39% of the original value.

*Acknowledgements*: The authors are grateful to two anonymous referees, whose comments helped to improve the exposition in this paper. This work started as a Group Project at the CBMS Conference *Interface of Mathematical Biology and Linear Algebra*. The authors extend their appreciation to the NSF and the CBMS, which sponsored that event. The authors are grateful to Zhuorong Mao and Poroshat Yazdanbakhshghahyazi for constructive conversations when this project was in its initial phase. Steve Kirkland's research is supported in part by the Natural Sciences and Engineering Research Council of Canada under grant number RGPIN–2019–05408.